# COMBINATORICS OF WORDS AND SEMIGROUP ALGEBRAS WHICH ARE SUMS OF LOCALLY NILPOTENT SUBALGEBRAS

Vesselin Drensky and Lakhdar Hammoudi

ABSTRACT. We construct new examples of non-nil algebras with any number of generators, which are direct sums of two locally nilpotent subalgebras. As all previously known examples, our examples are contracted semigroup algebras and the underlying semigroups are unions of locally nilpotent subsemigroups. In our constructions we make more transparent than in the past the close relationship between the considered problem and combinatorics of words.

1. INTRODUCTION

A fascinating area in ring theory is the study of rings that are sums of subrings with a given property. The interest in this area became obvious when Kegel proved that a ring which is a sum of two nilpotent subrings is itself nilpotent [K1]. The importance of studying such rings grew even bigger when Ferrero and Puczyłowski established that the famous Köthe conjecture [Koe] is equivalent to the statement that every ring which is a sum of a nilpotent subring and a nil subring must be nil [FP]. For some recent results on the Köthe conjecture and related topics we refer to the recent survey by Smoktunowicz [Sm]. Ferrero-Puczyłowski's work, as well as many others, are mainly motivated by possible generalizations of Kegel's result. Already, in [K2], Kegel asked whether a similar result to his could be drawn in the locally nilpotent and nil cases. It was only some 30 years later that Kelarev gave the first counterexample to Kegel's problem: A ring which is the sum of two locally nilpotent subrings maybe not nil [Ke1]. Nowadays, there are several other families of such rings [S1, S2, Ke2, F1, F2]. Although for all these examples every finitely generated ideal in each of the locally nilpotent summands is nilpotent, in reality, some of these rings may even contain a free subalgebra of rank 2 [F2]. All known examples are based on somewhat the same idea. The constructed ring is a contracted semigroup algebra $\mathbb{K}_0 S$ over a field $\mathbb{K}$, where the semigroup $S = \bar{S} \cup \{0\}$ is generated by two elements.

The first example of Kelarev [Ke1] is quite sophisticated. In the further examples, Kelarev [Ke2] and Fukshansky [F2] made use of a degree function on the free

Partially supported by Grant MM-1106/2001 of the Bulgarian Foundation for Scientific Research.
2000 *Mathematics Subject Classification.* 16N40, 16S15, 20M05, 20M25, 68R15

nonunitary semigroup $\langle x, y \rangle$. This function is defined inductively by $d(x) = 1$, $d(y) = -\alpha$ and $d(v_1 v_2) = d(v_1) + d(v_2)$ for all $v_1$, $v_2 \in \langle x, y \rangle$, where $\alpha > 0$ is a fixed irrational number. Then $S$ is the Rees factor semigroup of $\langle x, y \rangle$ modulo the ideal generated by all $z \in \langle x, y \rangle$ with $|d(z)| \geq C$ for some properly chosen constant $C$ (e.g. $C = 3$ for $\alpha = \sqrt{2}$ in [Ke2]).

Salwa's examples are based on semigroups of partial transformations of the real line [S1, S2]. It turns out that the considered semigroups are very close to semigroups defined by infinite words of a special kind. One can then use combinatorics of words to keep track on the growth and the Gelfand-Kirillov dimension of the corresponding algebras.

On the other hand, there is an old example of a two-generated nil semigroup which cannot be embedded into a nil algebra [Kol]. This semigroup as constructed by Kolotov is defined by an aperiodic infinite word and satisfies the identity $x^5 = 0$. The contracted semigroup algebra is of minimal possible growth and Gelfand-Kirillov dimension 2.

The purpose of the present note is to combine different ideas by Kolotov, Kelarev and Salwa to produce other examples of contracted semigroup algebras which are not nil but are direct sums of two locally nilpotent subalgebras. This provides yet other counterexamples to Kegel's problem. Since the underlying semigroups are Rees factors of free semigroups, our algebras are also monomial algebras, i.e. they have systems of defining relations consisting of monomials only. Our approach which is based on infinite words is quite natural and generalizes in a certain sense Kelarev's and Salwa's constructions. In the two-generated case our examples are "natural" homomorphic images of the examples of Kelarev and Fukshansky. Nevertheless, we think that they are interesting on their own and make more transparent the relationships with combinatorics of words. For example, we prove that any Sturmian word can be used to produce a non-nil ring which is the sum of two locally nilpotent subrings. We also show that the rings constructed by Salwa [S1, S2] are homomorphic images of Kelarev's [Ke2]. We shall start with some definitions and known facts.

## 2. Preliminaries

We denote by $\langle X \rangle$ the nonunitary free semigroup on a free set $X$ of generators. If $\mathbb{I}$ is an ideal of $\langle X \rangle$, then the Rees factor $\langle X \rangle / \mathbb{I}$ is the semigroup $(\langle X \rangle \setminus \mathbb{I}) \cup \{0\}$ with defining relations $w = 0$ if $w \in \mathbb{I}$. If $\Omega$ is any finite or infinite word on $X$, then we define a semigroup $S(\Omega)$ as the Rees factor $\langle X \rangle / \mathbb{I}$ where the ideal $\mathbb{I}$ coincides with the set of all words in $\langle X \rangle$ which are not subwords of $\Omega$. For any map $d : X \to \mathbb{R}$ we define a degree function on $\langle X \rangle$ (which we shall denote by the same symbol $d$) inductively by $d(w_1 w_2) = d(w_1) + d(w_2)$, $w_1, w_2 \in \langle X \rangle$. Clearly, $d$ induces a degree function on the nonzero elements of any Rees factor $\langle X \rangle / \mathbb{I}$ and again $d(uv) = d(u) + d(v)$ if $u, v \in \langle X \rangle / \mathbb{I}$ and $uv \neq 0$.

We fix a field $\mathbb{K}$ of any characteristics. If $S$ is a semigroup with zero, we denote by $\mathbb{K}_0 S$ the contracted semigroup algebra of $S$ over $\mathbb{K}$. This is the associative algebra which, as a vector space, has a basis $S \setminus \{0\}$ and the multiplication is induced by the multiplication of $S$. If $S = \langle X \rangle / \mathbb{I}$ for some ideal $\mathbb{I}$ of $\langle X \rangle$, then $\mathbb{K}_0 S$ is a monomial algebra. See the survey article by Belov, Borisenko and Latyshev [BBL] for different aspects of the theory of monomial algebras. For a subset $V$ of $\mathbb{K}_0 S$ we denote by $\mathrm{alg}_{\mathbb{K}_0 S}(V)$ and $\mathrm{id}_{\mathbb{K}_0 S}(V)$, respectively, the subalgebra and the ideal of $\mathbb{K}_0 S$ generated by $V$. Similarly, $\mathrm{semigr}(W)$ and $\mathrm{id}_S(W)$ are, respectively, the subsemigroup and the ideal of $S$ generated by the subset $W$ of $S$.

The following lemma appears (explicitly or implicitly) in [F1, F2, Ke1, Ke2, Kol, S1, S2] and is essential for the constructions there.

**Lemma 1.** (i) *Let $S$ be a semigroup with zero which is a union of its subsemigroups $S_1$ and $S_2$, $S_1 \cap S_2 = \{0\}$. Then, as a vector space, the algebra $\mathbb{K}_0 S$ is a direct sum of its subalgebras $\mathbb{K}_0 S_1$ and $\mathbb{K}_0 S_2$;*

(ii) *If the semigroup $S$ is locally nilpotent and its finitely generated ideals are nilpotent, then so is $\mathbb{K}_0 S$;*

(iii) *If $S = \langle X \rangle / \mathbb{I}$ for some finite set $X$ and some ideal $\mathbb{I}$ of $\langle X \rangle$ and $S$ is not nilpotent, then the algebra $\mathbb{K}_0 S$ is not nil.*

*Proof.* The statement in (i) is obvious and (ii) is routine: The elements of a finite subset $V = \{v_1, \ldots, v_k\}$ of $\mathbb{K}_0 S$ are linear combination of a finite subset $W = \{w_1, \ldots, w_n\}$ of $S$. Hence the subalgebra $\mathrm{alg}_{\mathbb{K}_0 S}(V)$ and the ideal $\mathrm{id}_{\mathbb{K}_0 S}(V)$ of $\mathbb{K}_0 S$ are contained, respectively, in the subalgebra and the ideal of $\mathbb{K}_0 S$ generated by $W$. Obviously, if $W \subset S$, then $\mathrm{alg}_{\mathbb{K}_0 S}(W) = \mathbb{K}_0 \mathrm{semigr}(W)$ and $\mathrm{id}_{\mathbb{K}_0 S}(W) = \mathbb{K}_0 \mathrm{id}_S(W)$. For (iii) it is sufficient to pay attention that if $X = \{x_1, \ldots, x_m\}$ and $x_{i_1} \cdots x_{i_n} \neq 0$ in $S$, then the element $(x_1 + \cdots + x_m)^n$ is a sum of all products $x_{j_1} \cdots x_{j_n}$ and contains the nonzero summand $x_{i_1} \cdots x_{i_n}$. Since $\mathbb{K}_0 S$ is a monomial algebra, we derive that $(x_1 + \cdots + x_m)^n \neq 0$.

## 3. Kelarev's and Salwa's examples and their generalizations

Our aim in this section is to give a "unified" construction to produce counterexamples to Kegel's problem by combining different ideas by Kelarev, Salwa and Fukshansky. Then we show that our method gives further counterexamples. For example, we prove that all semigroups of Sturmian words are union of two locally nilpotent subsemigroups. We shall start by recalling Kelarev's and Salwa's examples and proving that Salwa's rings are homomorphic images of the rings constructed as those by Kelarev in [Ke2]. We show that, in fact, they can be constructed using Kelarev's technique.

Let $\langle x, y \rangle$ be the nonunitary free semigroup over two variables. Fix an irrational number $\alpha > 0$ and a positive constant $C$. Consider the degree function $d : \langle x, y \rangle \to \mathbb{R}$ defined by $d(x) = 1$, $d(y) = -\alpha$ and $d(uv) = d(u) + d(v)$ for all $u, v \in \langle x, y \rangle$. Let $S$ be the factor semigroup $\langle x, y \rangle$ modulo the ideal generated by all words $z \in \langle x, y \rangle$ with $|d(z)| \geq C$. Define now the following subsemigroups:

$$S_1 = \{s \mid 0 \neq s \in S, d(s) > 0\} \cup \{0\}, \quad S_2 = \{s \mid 0 \neq s \in S, d(s) < 0\} \cup \{0\}.$$

In [Ke2] Kelarev proved that for any field $\mathbb{K}$ and for $C = 3$ and $\alpha = \sqrt{2}$, the semigroup algebra $\mathbb{K}_0 S$ is a non-nil ring which is the sum of the two locally nilpotent subrings $\mathbb{K}_0 S_1$ and $\mathbb{K}_0 S_2$. Moreover, Fukshansky proved that for properly chosen $\alpha$ and $C$ the semigroup algebras $\mathbb{K}_0 S$ contain a free algebra of rank 2 [F1, F2].

Salwa's examples are based on semigroups of translations on the real line. For any real $f \in \mathbb{R}$ and any open interval $I \subset \mathbb{R}$, denote by $(f, I)$ the translation $t \to t + f$, for all $t \in I$, where we assume that $(f, \emptyset) = 0$. For any two translations $(f, I)$ and $(g, J)$ their composition is the translation $(f, I) \circ (g, J) = (f + g, (I - g) \cap J)$. Let $b < 0 < a$ be two reals such that the fraction $a/b$ is irrational and fix an open interval $I \neq \emptyset$. Let $S$ be the semigroup of translations generated by $(a, I)$ and $(b, I)$ and consider the subsemigroups $S_1 = \{(f, J) \in S, J \neq \emptyset, f > 0\} \cup \{0\}$ and $S_2 = \{(f, J) \in S, J \neq \emptyset, f < 0\} \cup \{0\}$. In [S1] Salwa proved that for any interval $I$ such that $|I| \geq |a| + |b|$ and any field $\mathbb{K}$ of characteristics 0, the contracted semigroup algebra $\mathbb{K}_0 S$ is a non-nil algebra which is the sum of the two locally nilpotent subalgebras $\mathbb{K}_0 S_1$ and $\mathbb{K}_0 S_2$.

Now, let us define the degree function $d : \langle x, y \rangle \to \mathbb{R}$ by $d(x) = a$, $d(y) = b$ and $d(uv) = d(u) + d(v)$ for all $u, v \in \langle x, y \rangle$.

**Proposition 2.** *Let $S$ be the semigroup constructed by Salwa as above. Then $S$ is a homomorphic image of $\langle x, y \rangle / \mathbb{J}$, where $\mathbb{J}$ is the ideal of $\langle x, y \rangle$ generated by all words $xz, yz \in \langle x, y \rangle$ such that $|d(z)| > |I|$.*

*Proof.* Let $\tau$ be the homomorphism of semigroups defined by $\tau : \langle x, y \rangle \to S$, $\tau(x) = (a, I)$ and $\tau(y) = (b, I)$. We shall prove that $\mathbb{J} \subset \ker(\tau)$. Let $z = x_1 \cdots x_n \in \langle x, y \rangle$, where $x_i = x$ or $x_i = y$, $i = 1, \ldots, n$, be such that $|d(z)| > |I|$. Then $\tau(xz) = \tau(xx_1 \cdots x_n) = (a, I) \circ (f_1, I) \circ \cdots \circ (f_n, I) = (a + f_1 + \cdots + f_n, J)$, where $f_i = a$ or $f = b$, $i = 1, \ldots, n$, and $J$ is an interval such that $t \in J$ if and only if $t \in I$, $t + f_n \in I, \ldots, t + f_n + \cdots + f_1 \in I$, see [S1]. Therefore, $t \in J$ if and only if $t \in I$, $t \in I - f_n, \ldots, t \in I - f_n - \cdots - f_1$. Now, since $|d(z)| = |d(x_1) + \cdots + d(x_n)| = |f_1 + \cdots + f_n| > |I|$, we obtain that $(I - f_n - \cdots - f_1) \cap I = \emptyset$. So, there is no $t$ which is in $(I - f_n - \cdots - f_1) \cap I$. Hence $\tau(xz) = (a + f_1 + \cdots + f_n, \emptyset) = 0$. The same arguments work for elements of the form $yz$ with $|d(z)| > |I|$. So, $\mathbb{J} \subset \ker(\tau)$.

The following result combines the construction of Kelarev [Ke2] and Fukshansky [F2] with the ideas of Salwa [S1, S2].

**Theorem 3.** *Let $X$ and $Y$ be nonempty sets of any cardinality and let $d$ be a degree function on the free semigroup $F = \langle X, Y \rangle$ with the property that $d(x) > 0$, $d(y) < 0$ for all $x \in X$, $y \in Y$, and such that, if the finite sum*

$$\sum_{x_i \in X} p_i d(x_i) + \sum_{y_j \in Y} q_i d(y_i) = 0$$

*for some $p_i, q_j \in \mathbb{Q}$, $p_i, q_j \geq 0$, then $p_i = q_j = 0$ for all $x_i \in X$, $y_j \in Y$. Let $a, b \in \mathbb{R}$, $a, b > 0$, and let there exist $x_0 \in X$ and $y_0 \in Y$ such that $d(x_0) \leq a$, $d(y_0) \geq -b$. Then the semigroup*

$$S = F/\mathrm{id}_F \{w \mid |d(w)| \geq a + b\}$$

*is not nilpotent and is a union of the subsemigroups*

$$S_+ = \{s \in S \mid d(s) > 0\} \cup \{0\}, \ S_- = \{s \in S \mid d(s) < 0\} \cup \{0\}$$

*which are locally nilpotent, their finitely generated ideals are nilpotent and $S_+ \cap S_- = \{0\}$.*

*Proof.* Since the ideal $\mathrm{id}_F\{w \mid |d(w)| \geq a+b\}$ of $F$ contains all $w \in F$ with $|d(w)| \geq a+b$, the factor semigroup $S$ has no nonzero elements $s$ with $|d(s)| \geq a+b$. The properties of $d(x_i)$ and $d(y_j)$ guarantee that $S$ has no elements $s$ with $d(s) = 0$. Hence $S = S_+ \cup S_-$ and $S_+ \cap S_- = \{0\}$. If $s_1, \ldots, s_n$ is any finite set of nonzero elements in $S_+$, then there exists a positive $\varepsilon$ such that $d(s_i) > \varepsilon$, $i = 1, \ldots, n$. Hence $d(s_{i_1} \cdots s_{i_k}) > k\varepsilon$ and for $k \geq (a+b)/\varepsilon$ we have $s_{i_1} \cdots s_{i_k} = 0$. This implies that the subsemigroup of $S_+$ generated by $s_1, \ldots, s_n$ is nilpotent. Similarly, any nonzero element $s_0$ of the ideal of $S_+$ generated by $s_1, \ldots, s_n$ is of the form $s's_is''$ for some $s', s'' \in S_+$ and some $i = 1, \ldots, n$. Hence $d(s_0) = d(s') + d(s_i) + d(s'') > d(s_i) > \varepsilon$ and again the ideal $\mathrm{id}_{S_+}\{s_1, \ldots, s_n\}$ is nilpotent. The considerations for $S_-$ are analogous.

In order to see that $S$ is not nilpotent, it is sufficient to construct for any positive integer $n$ a product $w_n = z_1 \cdots z_n$ in $F$ of some not necessarily pairwise different variables $z_1, \ldots, z_n \in X \cup Y$ such that $w_n$ does not belong to the ideal of $F$ generated by all words $w$ with $|d(w)| \geq a+b$. We shall define $w_n$ inductively. Let $x_0 \in X$ and $y_0 \in Y$ be generators of $F$ with $0 < d(x_0) < a$ and $-b < d(y_0) < 0$. We choose $w_1 = z_1$, where $z_1 \in X \cup Y$ and $d(z_1)$ belongs to the open interval $(-b, a)$. For example, we may choose $w_1 = x_0$ or $w_1 = y_0$. If we have already fixed the element $w_{n-1}$ in $F$ with the property that $d(w_{n-1}) \in (-b, a)$, then there exists a $z_n \in X \cup Y$ such that $d(w_{n-1}) + d(z_n)$ also belongs to $(-b, a)$. For example, if $d(w_{n-1}) > 0$ we may choose $z_n = y_0$ and if $d(w_{n-1}) < 0$ we may fix $z_n = x_0$. Then we define $w_n = w_{n-1}z_n$ and again have $d(w_n) \in (-b, a)$. If $w_n \in \mathrm{id}_F\{w \in F \mid |d(w)| \geq a+b\}$, then $w_n = uwv$, where $u, v \in F \cup \{1\}$ and $|d(w)| \geq a+b$. Clearly $u = w_p$ and $uw = w_q$ for some $p < q$ (we may define formally $w_0 = 1$ and $d(w_0) = 0$, if $u = 1$). Hence $d(w) = d(uw) - d(u) = d(w_q) - d(w_p)$. Since $d(w_p), d(w_q) \in (-b, a)$, this immediately gives that $|d(w)| < a+b$ which is a contradiction. Hence the image of $w_n$ is not equal to 0 in $S$ and the semigroup $S$ is not nilpotent.

Let, in the notations and in the proof of Theorem 3, $w_n = z_1 \cdots z_n = w_{n-1}z_n \in F = \langle X, Y \rangle$, $n = 1, 2, \ldots$, be the constructed word of length $n$ which does not belong to the ideal $\mathrm{id}_F\{w \in F \mid |d(w)| \geq a+b\}$.

**Corollary 4.** *Let $\Omega = z_1 z_2 \cdots z_n \cdots$ be the infinite word in $X$ and $Y$ associated with the sequence of the words $w_n = w_{n-1}z_n$ from the proof of Theorem 3. Then the semigroup $S(\Omega)$ is a union of its subsemigroups*

$$S(\Omega)_+ = \{s \in S(\Omega) \mid d(s) > 0\} \cup \{0\}, \quad S(\Omega)_- = \{s \in S(\Omega) \mid d(s) < 0\} \cup \{0\}$$

*which are locally nilpotent, their finitely generated ideals are nilpotent and $S(\Omega)_+ \cap S(\Omega)_- = \{0\}$.*

*Proof.* The infinite word $\Omega$ does not contain finite subwords belonging to the ideal $\mathrm{id}_F\{w \in F \mid |d(w)| \geq a+b\}$. Hence $S(\Omega)$ is a homomorphic image of the semigroup $S = F/\mathrm{id}_F\{w \in F \mid |d(w)| \geq a+b\}$. Since $S(\Omega)$ is a Rees factor of $F$ modulo an ideal $\mathbb{I}$, we obtain that it inherits the degree function of $F$, $S(\Omega)_\pm$ is a homomorphic image of $S_\pm$ and $S(\Omega)_+ \cap S(\Omega)_- = \{0\}$. This gives that $S(\Omega)_\pm$ is locally nilpotent and its finitely generated ideals are nilpotent. Finally, the infinity of the word $\Omega$ gives that the semigroup $S(\Omega)$ is not nilpotent.

**Corollary 5.** *In the notations of Theorem 3 and Corollary 4, the contracted semigroup algebras $\mathbb{K}_0 S$ and $\mathbb{K}_0 S(\Omega)$ are not nil and are direct sums (as vector spaces) of their subalgebras $\mathbb{K}_0 S_+$, $\mathbb{K}_0 S_-$ and $\mathbb{K}_0 S(\Omega)_+$, $\mathbb{K}_0 S(\Omega)_-$, respectively. These four subalgebras are locally nilpotent and any of their finitely generated ideals is nilpotent.*

A Sturmian word is an infinite word which has exactly $n+1$ different subwords of length $n$ for all $n \geq 1$. By putting $n = 1$, one observes that any Sturmian word is written using an alphabet of two symbols. It is a well known combinatorial fact, see [CH] or [L, Chapter 1], that if an infinite word has less than $n_0 + 1$ subwords of length $n_0$ for some $n_0$, then it is eventually periodic. Then there exists a positive integer $c$ such that for all $n$ the word has less than $c$ different subwords of length $n$.

The proofs of the following claims about Sturmian words can be found in [L, Chapter 2]. Sturmian words can be nicely characterized using the notion of slope of a word. Let $w$ be a nonempty word in the variables $x$ and $y$. Then the slope of $w$ is the number $\pi(w) = \partial_x(w)/|w|$, where $|w|$ is the length of $w$ and $\partial_x(w)$ is the number of $x$'s appearing in $w$. The slope of the Sturmian word $\Omega$ is defined as the number $\alpha = \lim_{n\to\infty} \pi(u_n)$, where $(u_n)$ is the sequence of all nonempty initial subwords of $\Omega$. It is known that the slope of $\alpha$ is an irrational number in the interval $(0, 1)$. A strong connection between any nonempty subword $u$ of $\Omega$ and $\alpha$ is the inequality

$$|\pi(u) - \alpha| \leq \frac{1}{|u|}.$$

**Theorem 6.** *Over any field $\mathbb{K}$ and for any Sturmian word $\Omega$, the contracted semigroup algebra $\mathbb{K}_0 S(\Omega)$ is a non-nil algebra which is the sum of two locally nilpotent subalgebras $\mathbb{K}_0 S(\Omega)_+$ and $\mathbb{K}_0 S(\Omega)_-$.*

*Proof.* Let $\Omega$ be a Sturmian word with slope $\alpha$. We define on $S(\Omega)$ the degree function $d(u) = (1-\alpha)\partial_x(u) - \alpha\partial_y(u)$ for all $u \in S(\Omega)$. Consider now the subsemigroups

$$S(\Omega)_+ = \{s \mid 0 \neq s \in S,\ d(s) > 0\},\ S(\Omega)_- = \{s \mid 0 \neq s \in S,\ d(s) < 0\}.$$

The inequality $|\pi(u) - \alpha| \leq 1/|u|$ gives that

$$-1 \leq (1-\alpha)\partial_x(u) - \alpha\partial_y(u) \leq 1$$

and $|d(u)| \leq 1$ for all nonzero words $u \in S$. Clearly $d(x) = 1 - \alpha$, $d(y) = -\alpha$ and the fraction $d(x)/d(y)$ is irrational. Now we can complete the proof with the same arguments as in the proofs of Theorem 3 and Corollary 4.

It is important to observe that the algebras of Theorem 6 are extremal. For any Sturmian word $\Omega$, the growth function $g(n)$ of the contracted semigroup algebra $K_0 S(\Omega)$ (with respect to the generating set $\{x, y\}$) is $g(n) = n(n+3)/2$ and $\mathrm{GKdim}\,\mathbb{K}_0 S(\Omega) = 2$. The algebra $\mathbb{K}_0 S(\Omega)$ is a just-infinite dimensional monomial algebra: If we consider any proper homomorphic image $\overline{K_0 S(\Omega)}$ of $K_0 S(\Omega)$, then the corresponding growth function $\bar{g}(n)$ satisfies $\bar{g}(n_0) < n_0(n_0+3)/2$ for some $n_0$ and by the Bergman gap theorem [B], see also [KL], the algebra $\overline{K_0 S(\Omega)}$ would be nilpotent. One can also observe that the Jacobson radical $JK_0 S(\Omega)$ of $K_0 S(\Omega)$ is zero. Indeed, if we assume that $JK_0 S(\Omega) \neq \langle 0 \rangle$, then, since $K_0 S(\Omega)$ is a monomial algebra, $JK_0 S(\Omega)$ is locally nilpotent [BF]. But $K_0 S(\Omega)$ is just-infinite. Thus, $JK_0 S(\Omega)$ is of finite codimension. It is obvious now, that $K_0 S(\Omega)/(JK_0 S(\Omega))$ is a finite dimensional algebra in which every monomial is nilpotent. So, it is nilpotent and therefore, $K_0 S(\Omega)$ is also nilpotent. This contradicts the construction of $K_0 S(\Omega)$ and proves our claim.

## 4. Kolotov's example and some more generalizations

Salwa [S1, S2] showed that some of his examples have homomorphic images of the form $S(\Omega)$ where the infinite word $\Omega$ is Sturmian and has some extremal properties (cf. [2, Chapter 2]). In Section 3 we proved that the semigroup of any Sturmian word is the union of locally nilpotent subsemigroups (Theorem 6). Here we shall study concrete examples of Sturmian words arising from Kolotov's example [Kol]. This class of words satisfies the hypothesis of Corollary 4. We then look at some generalizations such as semigroups with arbitrary number of generators.

Slightly changing the notations, Kolotov's semigroup is defined in the following way: Let $X = \{x, y\}$ and
$$\omega_0 = x, \ \omega_1 = xy, \ \omega_{n+2} = \omega_{n+1}\omega_n\omega_{n+1}, \ n = 0, 1, 2, \ldots$$
Then $\Omega$ is the infinite word obtained as a limit (to the right) of the words $\omega_n$ (i.e. each $\omega_n$, $n = 0, 1, 2, \ldots$, is a beginning of $\Omega$). Kolotov proved that the semigroup $S(\Omega)$ satisfies the identity $w^5 = 0$ and contains exactly $n+1$ different nonzero elements of length $n$. Thus, it is a Sturmian word.

Now we generalize Kolotov's semigroup to an infinite family of semigroups with an arbitrary finite number of generators. Let $\langle Z_m \rangle = \langle z_1, \ldots, z_m \rangle$ be the nonunitary free semigroup with $m > 1$ generators. For any word $w = z_{i_1} \cdots z_{i_n}$, we denote by
$$\bar{w} = \bar{z}_1^{c_1} \cdots \bar{z}_m^{c_m}$$
its abelianization (i.e. $\bar{w}$ is $w$ modulo the commutativity relations $z_i z_j = z_j z_i$, $i, j = 1, \ldots, m$). We associate with $w$ the polynomial with nonnegative integer coefficients
$$f_w(t) = c_m t^{m-1} + c_{m-1} t^{m-2} + \cdots + c_2 t + c_1 \in \mathbb{Z}[t].$$

We fix an infinite sequence of words $\rho_n = \rho_n(z_1, \ldots, z_m) \in \langle Z_m \rangle$, $n = 0, 1, 2, \ldots$ with the same abelianizations $\bar{\rho} = \bar{\rho}_n = \bar{z}_1^{b_1} \cdots \bar{z}_m^{b_m}$, $n \geq 0$, and consider arbitrary (but fixed) $m$ words $\omega_0, \omega_1, \ldots, \omega_{m-1} \in \langle Z_m \rangle$ with abelianizations $\bar{\omega}_j = \bar{z}_1^{a_{1j}} \cdots \bar{z}_m^{a_{mj}}$, $j = 0, 1, \ldots, m-1$. We define the words $\omega_n$ by the "linear" recurrent relations
$$\omega_{n+m} = \rho_n(\omega_n, \omega_{n+1}, \ldots, \omega_{n+m-1}), \ n = 0, 1, 2, \ldots.$$
Consider now the ideal $I(\omega, \rho)$ of $\langle Z_m \rangle$ consisting of all words $w$ which are not subwords of any word $\omega_n$, $n = 0, 1, 2, \ldots$. We denote by $S(\omega, \rho)$ the Rees factor semigroup $\langle Z_m \rangle / I(\omega, \rho)$.

**Theorem 7.** *Let $m > 1$ and let the sequence of words $\rho_n$, $n = 0, 1, 2, \ldots$, with the same abelianization $\bar{\rho}$ be chosen in such a way that the associated polynomial $t^m - f_\rho(t) \in \mathbb{Z}[t]$ is irreducible over $\mathbb{Q}$ and has a real zero in the interval $(-1, 0)$. Suppose that the $m \times m$ determinant $\det(a_{ij})$ associated with the abelianization of the words $\omega_0, \omega_1, \ldots, \omega_{m-1}, \rho_i \in \langle Z_m \rangle$ is different from zero. Then the semigroup $S(\omega, \rho)$ is not nilpotent and there exists a degree function on $\langle Z_m \rangle$ such that $S(\omega, \rho)$ is a union of the subsemigroups $S(\omega, \rho)_+$ and $S(\omega, \rho)_-$ which are locally nilpotent, their finitely generated ideals are nilpotent and $S(\omega, \rho)_+ \cap S(\omega, \rho)_- = \{0\}$.*

*Proof.* We shall prove the theorem in several steps.

*Step 1.* Let $\xi \in (-1, 0)$ be a zero of the polynomial $t^m - f_\rho(t)$. Since $t^m - f_\rho(t)$ is irreducible over $\mathbb{Q}$, the set $\{1, \xi, \xi^2, \ldots, \xi^{m-1}\}$ is a basis of $\mathbb{Q}(\xi) \subset \mathbb{R}$. We are looking for a degree function $d : \langle Z_m \rangle \to \mathbb{Q}(\xi)$ such that $d(\omega_j) = \xi^j$, $j = 0, 1, \ldots, m-1$. We consider the equalities

$$d(\omega_j) = a_{1j}d(z_1) + a_{2j}d(z_2) + \cdots + a_{mj}d(z_m) = \xi^j, \ j = 0, 1, \ldots, m-1,$$

as a system of $m$ linear equations with $m$ unknowns $d(z_1), d(z_2), \ldots, d(z_m)$. Since the determinant $\det(a_{ij})$ is different from zero, the system has a unique solution in $\mathbb{Q}(\xi)$. The linear combinations $d(\omega_j) = \sum_{i=0}^{m-1} a_{ij}d(z_i) = \xi^j$ of $d(z_1), d(z_2), \ldots, d(z_m)$, $j = 0, 1, \ldots, m-1$, are linearly independent over $\mathbb{Q}$. Hence $d(z_1), d(z_2), \ldots, d(z_m)$ are also linearly independent over $\mathbb{Q}$ and this implies that $\langle Z_m \rangle$ contains no elements $w$ with $d(w) = 0$.

*Step 2.* We shall show that $d(\omega_{n+m}) = \xi^{n+m}$, $n = 0, 1, 2, \ldots$. Since $\bar{\rho}(\bar{z}_1, \ldots, \bar{z}_m) = \bar{z}_1^{b_1} \cdots \bar{z}_m^{b_m}$, and by induction, bearing in mind that $\xi^m - f_\rho(\xi) = 0$, we obtain that

$$d(\omega_{n+m}) = b_m d(\omega_{n+m-1}) + b_{m-1} d(\omega_{n+m-2}) + \cdots + b_2 d(\omega_{n+1}) + b_1 d(\omega_n) = \xi^n f_\rho(\xi) = \xi^{n+m}.$$

*Step 3.* Let $c = b_1 + b_2 + \cdots + b_m$ be the length of the words $\rho_n(z_1, \ldots, z_m) \in \langle Z_m \rangle$ and let all beginnings of the words $\omega_0, \omega_1, \ldots, \omega_n$ have degrees belonging to the interval $[-q_n, q_n]$ for some positive number $q_n$. We shall show that $q_{n+m} \leq (c-1)|\xi|^n + q_{n+m-1}$ for all $n = 0, 1, 2, \ldots$. For this purpose, if $\rho_n(z_1, \ldots, z_m) = z_{i_1} \cdots z_{i_c}$, then $\omega_{n+m} = \omega_{n-1+i_1} \cdots \omega_{n-1+i_c}$. Hence any beginning $u_{n+m}$ of $\omega_{n+m}$ is of the form $u_{n+m} = \omega_{n-1+i_1} \cdots \omega_{n-1+i_k} u_{n+m-1}$, where $k = 0, 1, \ldots, c-1$ and $u_{n+m-1}$ is a beginning of some word $\omega_0, \omega_1, \ldots, \omega_{n+m-1}$. Thus,

$$d(u_{n+m}) = d(\omega_{n-1+i_1}) + \cdots + d(\omega_{n-1+i_k}) + d(u_{n+m-1}),$$

$$|d(u_{n+m})| \leq |d(\omega_{n-1+i_1})| + \cdots + |d(\omega_{n-1+i_k})| + |d(u_{n+m-1})| \leq$$

$$\leq k|\xi|^n + q_{n+m-1} \leq (c-1)|\xi|^n + q_{n+m-1}.$$

*Step 4.* Starting from $n = 0$, by induction we obtain that we may take

$$q_{n+m} = (c-1)(1 + |\xi| + |\xi|^2 + \cdots + |\xi|^n) + q_{m-1} < \frac{c-1}{1-|\xi|} + q_{m-1} = q.$$

So, the absolute values of the degrees of the beginnings of all words $\omega_n$, $n = 0, 1, 2, \ldots$, are bounded by the same bound $q$. We can decompose $S(\omega, \rho)$ as the union of the locally nilpotent subsemigroups (which are with nilpotent finitely generated ideals) $S(\omega, \rho)_+$ and $S(\omega, \rho)_-$ with $S(\omega, \rho)_+ \cap S(\omega, \rho)_- = \{0\}$.

*Step 5.* In order to show that the semigroup $S(\omega, \rho)$ is not nilpotent, it is sufficient to see that there are words $\omega_n$ of arbitrarily large length. Since $|\xi| < 1$, the numbers $\xi^n$, $n = 0, 1, 2, \ldots$, are pairwise different. The words $\omega_n$ are also pairwise different, because $d(\omega_n) = \xi^n$. There are only $m^k$ words of length $k$ in $\langle Z_m \rangle$. Hence there are words $\omega_n$ of arbitrarily large length.

**Examples 8.** (i) Kolotov's semigroup can be obtained from Theorem 7 for $m = 2$, $\rho_n = \rho = z_2 z_1 z_2$, $n = 0, 1, 2, \ldots$, $\omega_0 = z_1$ and $\omega_1 = z_1 z_2$. The determinant is equal to
$$\Delta = \begin{vmatrix} 1 & 0 \\ 1 & 1 \end{vmatrix}.$$

The polynomial $t^2 - (2t + 1)$ is irreducible over $\mathbb{Q}$ and has a zero $1 - \sqrt{2} \in (-1, 0)$.

(ii) In order to produce other examples for any $m > 1$, we need a set of words $\{\omega_0, \omega_1, \ldots, \omega_{m-1}\}$ with nonzero associated determinant (easy to arrange) and a polynomial $f_\rho(t)$ with nonnegative integer coefficients such that $t^m - f_\rho(t)$ is irreducible over $\mathbb{Q}$ and $(-1)^m - f_\rho(-1) > 0$ (obviously $0^m - f_\rho(0) < 0$). Such a polynomial may be arranged in the following way. We fix a prime $p$ and consider a polynomial $f_\rho(t) = p(c_m t^{m-1} + c_{m-1} t^{m-2} + \cdots + c_2 t + c_1)$ where $c_1, c_2, \ldots, c_m$ are positive integers and $p$ does not divide $c_1$. By the Eisenstein criterion the polynomial $t^m - f_\rho(t)$ is irreducible over $\mathbb{Q}$. Choosing $c_2$ sufficiently large, we can guarantee that $(-1)^m - f_\rho(-1) > 0$.

It is worth noting that different examples of monomial algebras with arbitrary Gelfand-Kirillov dimension and based on infinite words defined by linear recurrent relations were constructed by Vishne [V].

Institute of Mathematics and Informatics, Bulgarian Academy of Sciences, Acad. G. Bonchev Str., Block 8, 1113 Sofia, Bulgaria
*E-mail address*: `drensky@math.bas.bg`

Department of Mathematics, Ohio University, 571 West Fifth Street, Chillicothe, OH 45601 USA
*E-mail address*: `hammoudi@ohio.edu`